\input amstex
\documentstyle{amsppt}
\nologo
\NoBlackBoxes
%%\NoRunningHeads
%\mag=1200
%\hsize=31 pc  
%\vsize=44 pc  
%\hcorrection{5mm}
%\baselineskip 15pt

\hcorrection{19mm}

\topmatter
\title    Annular Dehn fillings
\endtitle
\author    Cameron McA. Gordon$^{1}$ and Ying-Qing Wu$^2$
\endauthor
\leftheadtext{C.~M{\fiverm c}A.~Gordon and Y-Q.~Wu}
\rightheadtext{Annular Dehn fillings}
\address Department of Mathematics, University of Texas at Austin,
Austin, TX 78712
\endaddress
\email  gordon\@math.utexas.edu
\endemail
\address Department of Mathematics, University of Iowa, Iowa City, IA
52242
\endaddress
\email  wu\@math.uiowa.edu
\endemail
\thanks  $^1$ Partially supported by NSF grant \#DMS 9626550.
\endthanks 
\thanks  $^2$ Partially supported by NSF grant \#DMS 9802558.
\endthanks

\abstract
We show that if a simple 3-manifold $M$ has two Dehn fillings at
distance $\Delta \geq 4$, each of which contains an essential annulus,
then $M$ is one of three specific 2-component link exteriors in
$S^3$.  One of these has such a pair of annular fillings with $\Delta
= 5$, and the other two have pairs with $\Delta = 4$.  
\endabstract

\endtopmatter
 
\document
\define\proof{\demo{Proof}}
\define\endproof{\qed \enddemo}

\define\a{\alpha}
\redefine\b{\beta}
\redefine\tilde{\widetilde}
\redefine\hat{\widehat}

\redefine\bdd{\partial}

\define\Int{\text{\rm Int}}
\define\val{\text{\rm val}}

\input epsf.tex

\head \S 1.  Introduction
\endhead

Let $M$ be a (compact, connected, orientable) 3-manifold with a torus
boundary component $T_0$.  If $r$ is a {\it slope\/} (the isotopy
class of an essential unoriented simple loop) on $T_0$, then as usual
we denote by $M(r)$ the 3-manifold obtained from $M$ by {\it $r$-Dehn
filling}, that is, attaching a solid torus $J$ to $M$ along $T_0$ in
such a way that $r$ bounds a meridian disk in $J$.  

We shall say that a compact, connected, orientable $3$-manifold $M$ is
{\it simple\/} if it contains no essential surface of non-negative
Euler characteristic, i.e., sphere, disk, annulus or torus.  If $M$
has non-empty boundary and is not the 3-ball, then $M$ is simple if
and only if $M$ with its boundary tori removed has a complete
hyperbolic structure of finite volume with totally geodesic boundary
[Th1, Th2].  If $M$ is closed, then the geometrization conjecture
asserts that $M$ is simple if and only if $M$ is either hyperbolic or
belongs to a certain small class of Seifert fiber spaces [Th1, Th2].

If $M$ is hyperbolic, then Dehn fillings on $M$ are hyperbolic if we
exclude finitely many slopes from each torus boundary component [Th1,
Th2].  By doubling $M$ along its non-torus boundary components, we see
that if $M$ is simple then $M(r)$ is simple for all but finitely many
slopes $r$ on any given torus boundary component $T_0$, and a good
deal of attention has been directed towards obtaining a more precise
quantification of this statement.  Denote by $\Delta (r_1, r_2)$ the
distance, or minimal geometric intersection number, between two slopes
$r_1, r_2$ on a torus.  Define a 3-manifold to be of type $S$, $D$,
$A$ or $T$ if it contains an essential sphere, disk, annulus or torus,
respectively.  For $X_i \in \{S, D, A, T\}$, $i=1,2$, define
$\Delta(X_1, X_2)$ to be the maximum of $\Delta(r_1, r_2)$, where
$r_1$ and $r_2$ are Dehn filling slopes of some simple manifold $M$
such that $M(r_i)$ is of type $X_i$.  These numbers $\Delta (X_1,
X_2)$ are now known in all ten cases; see [GW2] for more details.

Except when $(X_1, X_2) = (A,A)$, $(A,T)$ or $(T,T)$, it is also
known that $\Delta(X_1, X_2)$ is realized by infinitely many
simple manifolds $M$; see [EW].  On the other hand, $\Delta(T, T)
= 8$, and there are exactly two simple manifolds $M$ admitting
toroidal fillings $M(r_1)$, $M(r_2)$ with $\Delta = \Delta(r_1, r_2) =
8$, exactly one with $\Delta = 7$, exactly one with $\Delta = 6$, and
infinitely many with $\Delta = 5$ [Go].  Similarly, $\Delta(A, T) =
5$ [Go, GW1], and there is exactly one simple manifold $M$ having
an annular filling $M(r_1)$ and a toroidal filling $M(r_2)$ with
$\Delta = 5$, exactly two with $\Delta = 4$, and infinitely many with
$\Delta = 3$ [GW1].  In the present paper we complete the picture by
dealing with the case $(A, A)$.  In this case, $\Delta(A, A) = 5$ [Go,
GW1], and there are infinitely many simple manifolds $M$ admitting
annular fillings $M(r_1), M(r_2)$ with $\Delta = 3$ [GW1].  Here we
show that there is exactly one such manifold $M$ with $\Delta = 5$,
and exactly two with $\Delta = 4$.  More precisely, we have the
following theorem.

\proclaim{Theorem 1.1} Suppose $M$ is a compact, connected,
orientable, irreducible, $\bdd$-irreducible, anannular $3$-manifold
which admits two annular Dehn fillings $M(r_1)$, $M(r_2)$ with $\Delta
= \Delta(r_1, r_2) \geq 4$.  Then one of the following holds.

(1)  $M$ is the exterior of the Whitehead link, and $\Delta = 4$.

(2)  $M$ is the exterior of the $2$-bridge link associated to the
rational number $3/10$, and $\Delta = 4$.

(3)  $M$ is the exterior of the $(-2,3,8)$ pretzel link, and $\Delta =
5$.  
\endproclaim

The three manifolds listed in the theorem are the exteriors of the
links in $S^3$ shown in Figure 1.1.

\bigskip
\leavevmode

\epsfxsize=5in      % adjust the width  of the figure
\centerline{\epsfbox{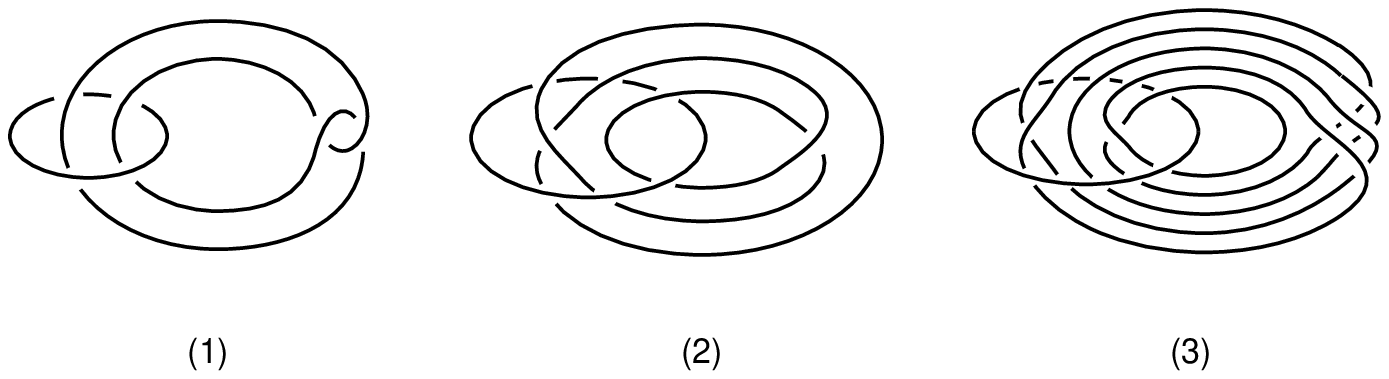}} 
\bigskip
\centerline{Figure 1.1}
\bigskip

That each of these link exteriors does have a pair of annular fillings
with $\Delta = 4$, $4$ and $5$ respectively is proved in [GW1].  The
fillings in question are also toroidal [GW1], so in fact these are
exactly the same manifolds which admit an annular and a toroidal
filling with $\Delta \geq 4$ [GW1, Theorem 1.1].  Using [GW1, Theorem
1.1], Qiu has independently proved Theorem 1.1 in the special case
where $M$ is the exterior of a knot in a solid torus [Q].

According to the proof of [GW1, Theorem 7.5], the annular fillings on
the three manifolds listed in Theorem 1.1 are non Seifert fibered
graph manifolds.  If $M$ admits some Seifert fibered surgery, then
$\bdd M$ consists of tori, in which case $M$ is hyperbolic if and only
if it is simple.  Hence the following corollary is an immediate
consequence of Theorem 1.1.

\proclaim{Corollary 1.2} Suppose $M$ is a compact orientable
hyperbolic 3-manifold with at least two torus boundary components, and
suppose $M(r_1), M(r_2)$ are Seifert fibered manifolds.  Then $\Delta
(r_1, r_2) \leq 3$.  \endproclaim

The condition that $M$ has at least two boundary components cannot be
removed.  For example, if $M$ is the figure 8 knot complement, then
$M(3)$ and $M(-3)$ are Seifert fibered, and $\Delta (-3, 3) = 6$.  It
is not known whether the bound 3 in the corollary is the best possible.

The proof of Theorem 1.1 proceeds as follows.  For $\alpha = 1,2$, let
$A_{\a}$ be an essential annulus in $M(r_{\a})$, meeting the Dehn
filling solid torus $J_{\a}$ in $n_{\a}$ meridian disks, with $n_{\a}$
minimal over all choices of $A_{\a}$.  This gives rise to a punctured
annulus $F_{\a} = A_{\a} \cap M$ in $M$, such that the boundary
components of $F_{\a}$ which lie on $T_0$ have slope $r_{\a}$, $\a =
1,2$.  The arcs of intersection of $F_1$ and $F_2$ then define labeled
graphs $G_{\a}$ in $A_{\a}$ with $n_{\a}$ vertices, $\a = 1,2$.  We
assume that $\Delta = \Delta(r_1, r_2) = 4$ or $5$, and show by a
detailed analysis that there are only three such pairs of graphs,
corresponding to the three examples listed in the theorem.  

The paper is organized as follows.  In Section 2 we give some
definitions and establish some basic properties of the graphs
$G_{\a}$.  In Section 3 we show that any graph in an annulus with no
trivial loops or parallel edges must satisfy one of four
possibilities; if the reduced graph $\hat G_{\a}$ is of the fourth
type we say that $G_{\a}$ is {\it special}.  Section 4 is devoted to
showing that if one of the graphs $G_1, G_2$ is special then they both
are, and (up to relabeling) $n_1 = 1$, $n_2 = 2$.  Section 5 considers
the generic case, $n_1, n_2 > 2$.  This is shown to be impossible, by
eliminating in turn the first three possibilities of Section 3 for the
reduced graphs $\hat G_{\a}$.  Section 6 shows that the case $n_1 =
2$, $n_2 > 2$ is also impossible.  In Section 7 we show that if $G_1$
and $G_2$ are special, so $n_1 = 1$ and $n_2 = 2$, then there is
exactly one possible pair $G_1, G_2$, with $\Delta = 4$, corresponding
to case (1) of Theorem 1.1.  Finally in Section 8, we show that if
$G_1, G_2$ are not special and $n_1, n_2 \leq 2$, then there are
exactly two possible pairs $G_1, G_2$, one with $\Delta = 4$, $n_1 =
n_2 = 2$, and one with $\Delta = 5$, $n_1 = n_2 = 2$, corresponding to
cases (2) and (3) of Theorem 1.1.

We would like to thank the referee for his/her careful reading and
helpful comments.

\head \S 2.  Preliminary Lemmas
\endhead

Throughout this paper, we will always assume that $M$ is a compact,
connected, irreducible, $\bdd$-irreducible, anannular 3-manifold, with
a torus boundary component $T_0$.  We use $\a, \b$ to denote the
numbers $1$ or $2$, with the convention that if they both appear, then
$\{\a, \b\} = \{1,2\}$.  Let $r_1, r_2$ be slopes on $T_0$ such that
$M(r_1), M(r_2)$ are annular, and let $A_{\a}$ be an essential annulus
in $M(r_{\a})$ such that $n_{\a}$, the number of components of
intersection of $A_{\a}$ with the Dehn filling solid torus $J_{\a}$,
is minimal among all essential annuli in $M(r_{\a})$, $\a = 1,2$.
Denote by $F_{\a}$ the punctured annulus $A_{\a } \cap M$.  Denote by
$\Delta = \Delta(r_1, r_2)$ the minimal geometric intersection number
between $r_1$ and $r_2$.  By [Go, Theorem 1.3] we have $\Delta \leq
5$.  Throughout this paper, we will always assume $\Delta = 4$ or $5$,
unless otherwise stated.

Minimizing the number of components of $F_1 \cap F_2$ by an isotopy,
we may assume that $F_1\cap F_2$ consists of arcs and circles which
are essential on both $F_{\a}$.  Let $u_1,\ldots, u_{n_{\a}}$ be the
disks that are the components of $A_{\a}\cap J_{\a}$, labeled
successively when traveling along $J_{\a}$.  Similarly let $v_1,
\ldots, v_{n_{\b}}$ be the disks in $A_{\b}\cap J_{\b}$.  Let $G_{\a}$
be the graph on $A_{\a}$ with the $u_i$'s as (fat) vertices, and the
arc components of $F_1\cap F_2$ with at least one endpoint on $T_0$ as
edges.  Note that we do not regard an edge endpoint on the boundary of
the annulus as a vertex, so we are abusing terminology somewhat in
that our graphs may have edge endpoints that do not lie on vertices.
The minimality of the number of components in $F_1\cap F_2$ and the
minimality of $n_{\a}$ imply that $G_{\a}$ has no trivial loops, and
that each disk face of $G_{\a}$ in $A_{\a}$ has interior disjoint from
$F_{\b}$.

An edge $e$ of a graph $G$ on an annulus $A$ is a {\it boundary
edge\/} if it has one endpoint on the boundary of $A$, otherwise it is
an {\it interior edge}.  A vertex $v$ of $G$ is a {\it boundary vertex
\/} if it is incident to a boundary edge, otherwise it is an {\it
interior vertex}.  Similarly, a face of $G$ is a {\it boundary face\/}
if it contains a boundary edge.

If $e$ is an edge of $G_{\a}$ with an endpoint $x$ on a fat vertex
$u_i$, then $x$ is labeled $j$ if $x$ is in $\bdd u_i\cap \bdd v_j$.
When going around the boundary of a vertex in
$G_{\a}$, the labels of the edge endpoints appear as $1, 2, \ldots,
n_{\b}$ repeated $\Delta$ times.

An edge $e$ at a vertex $u_i$ of $G_{\a}$ is called a {\it $j$-edge\/}
at $u_i$ if it has an endpoint at $u_i$ labeled $j$.  Dually, a
$j$-edge at $u_i$ is also an $i$-edge at $v_j$ in $G_{\b}$.  We say
that $e$ is an {\it $(i,k)$-edge\/} if it has labels $i$ and $k$ at
its two endpoints.

Each vertex of $G_{\a}$ is given a sign according to whether $J_{\a}$
passes $A_{\a}$ from the positive side or negative side at this
vertex.  Two vertices of $G_{\a}$ are {\it parallel\/} if they have
the same sign, otherwise they are {\it antiparallel.}  Note that if
$A_{\a}$ is a separating surface in $M(r_{\a})$, then $n_{\a}$ is
even, and $u_i, u_j$ are parallel if and only if $i, j$ have the same
parity.  An interior edge of $G_{\a}$ is a {\it positive edge\/} if it
connects parallel vertices.  Otherwise it is a {\it negative edge}.
We use $\val(v, G)$ to denote the valency of a vertex $v$ in a graph
$G$.

By considering each family of parallel edges of $G_{\a}$ as a single
edge $E$, we get the {\it reduced graph\/} $\hat{G}_{\a}$ on $A_{\a}$.
It has the same vertices as $G_{\a}$.  Denote by $|E|$ the number of
edges in $G_{\a}$ represented by $E$.

A cycle in $G_{\a}$ consisting of positive edges is a {\it Scharlemann
cycle\/} if it bounds a disk with interior disjoint from the graph,
and all the edges in the cycle have the same pair of labels $(i, i+1)$
at their two endpoints.  ($i+1 = 1$ if $i=n_{\b}$.)  The pair $(i,
i+1)$ is called the {\it label pair\/} of the Scharlemann cycle.  In
particular, a pair of adjacent parallel positive edges with the same
label pair is a Scharlemann cycle.  The boundary of the disk $D$
bounded by a Scharlemann cycle consists of edges of the Scharlemann
cycle and some arcs on the annulus $C_i$ on $T_0$ between $\bdd v_i$
and $\bdd v_{i+1}$.  When $n_{\b} = 2$, the two annuli $C_1$ and $C_2$
are still distinct, allowing one to differentiate between a
$(1,2)$-Scharlemann cycle and a $(2,1)$-Scharlemann cycle.  A pair of
edges $\{e_1, e_2\}$ is an {\it extended Scharlemann cycle\/} if there
is a Scharlemann cycle $\{e'_1, e'_2\}$ such that $e_i$ is parallel
and adjacent to $e'_i$.

A subgraph $G'$ of a graph $G$ on a surface $F$ is {\it essential\/}
if it is not contained in a disk in $F$. 

\proclaim{Lemma 2.1} 
(1) {\rm (The Parity Rule)} An edge $e$ is a positive edge in $G_1$ if
and only if it is a negative edge in $G_2$.

(2) A pair of edges cannot be parallel on both $G_1$ and $G_2$.

(3) If $G_{\a}$ has a set of $n_{\b}$ parallel negative edges, then on
$G_{\b}$ they form mutually disjoint essential cycles of equal length.

(4) If $G_{\a}$ has a Scharlemann cycle, then $A_{\b}$ is
separating, and $n_{\b}$ is even.  Moreover, the edges of the
Scharlemann cycle and the vertices at their endpoints form an
essential subgraph of $G_{\b}$.

(5)  $G_{\a}$ contains no extended Scharlemann cycle.
\endproclaim

\proof See [GW1, Lemma 2.2], except for (2) in the case that the pair
of edges $e_1, e_2$ are boundary edges.  If $e_1, e_2$ are boundary
edges parallel on both $G_1, G_2$, then they cut off bands $B_1, B_2$
on the punctured annuli $F_1, F_2$, which can be glued together to get
an annulus in the manifold $M$, which intersects the Dehn filling
torus $T_0$ in an essential circle.  This contradicts the assumption
that $M$ is $\bdd$-irreducible and anannular.  \endproof

Let $E$ be an edge of $\hat G_{\a}$ representing $n_{\b}$ parallel
negative edges on $G_{\a}$, connecting $u_i$ to $u_j$.  Then $E$
defines a permutation $\varphi: \{1, \ldots, n_{\b}\} \to \{1, \ldots,
n_{\b}\}$, such that an edge $e$ in $E$ has label $k$ at $u_i$ if and
only if it has label $\varphi(k)$ at $u_j$.  Call $\varphi$ the {\it
permutation associated to $E$}.  Because of the ambiguity in the order
of $u_i, u_j$, the permutation is only well defined up to inverse.  An
{\it $E$-orbit\/} is an orbit of $\varphi$.  Such an orbit determines
a cycle in $G_{\b}$ consisting of the edges of $E$ with endpoint
labels in this orbit, called the cycle of this orbit.  Note that all
the vertices in a cycle are parallel.  Topologically each such cycle
is a circle.  Lemma 2.1(3) says that these circles are mutually
disjoint, mutually parallel, essential circles on the annulus
$A_{\b}$.

\proclaim{Lemma 2.2} 
(1) Any two Scharlemann cycles on $G_{\a}$ have the same label pair.

(2) If $E$ is a positive edge in $\hat G_{\a}$, then $|E| \leq
n_{\b}/2 + 1$.  Moreover, if $|E| = n_{\b}/2 + 1$, then the
corresponding edges of $G_{\a}$ contain a Scharlemann cycle.

(3) Any family of parallel interior edges in $G_{\a}$ contains at
most $n_{\b}$ edges.  \endproclaim

\proof 
See [GW1, Lemma 2.5].
\endproof

\proclaim{Lemma 2.3}
(1)  If some vertex of $G_{\a}$ has more than $n_{\b}$ negative
edge endpoints, then $G_{\b}$ contains a Scharlemann cycle.

(2)  No vertex of $G_{\a}$ has more than $2n_{\b}$ negative edge
endpoints.  
\endproclaim

\proof 
For any label $i$ of $G_{\b}$, let $G_{\b}^+(i)$ be the subgraph of
$G_{\b}$ consisting of all vertices of $G_{\b}$ and all positive
$i$-edges of $G_{\b}$.  The edges of $G_{\b}^+(i)$ correspond to
the negative edges of $G_{\a}$ incident to the vertex $u_i$.  Let
the number of such edges be $k$.  Then if $f$ denotes the sum of the
Euler characteristics of the faces of $G_{\b}^+(i)$, we have 
$$ 0 = \chi(A_{\b}) = n_{\b} - k + f.$$
Therefore, if $k > n_{\b}$, $G_{\b}^+(i)$ has a disk face $D$.  Then
there is a Scharlemann cycle of $G_{\b}$ in $D$ by [HM, Proposition
5.1].  This proves (1).

To prove (2), assume $k > 2n_{\b}$.  Then by the above we have $f = k
- n_{\b} > n_{\b}$, so $G_{\b}^+(i)$ has more than $n_{\b}$ disk
faces, and by [HM, Proposition 5.1] each such face contains a
Scharlemann cycle of $G_{\b}$.  Hence $G_{\b}$ contains $s > n_{\b}$
Scharlemann cycles, all on the same label pair, say $(1,2)$, by
Lemma 2.2(1).  Define a graph $H$ in $A_{\b}$ as follows; see [GL,
Proof of Theorem 2.3].  The vertices of $H$ consist of the vertices of
$G_{\b}$, together with a vertex $v_{D}$ in the interior of each disk
face of $G_{\b}$ bounded by a Scharlemann cycle.  The edges of $H$ are
defined by joining each vertex $v_{D}$, within $D$, to the vertices of
$G_{\b}$ in $\bdd D$.  Thus $H$ has $n_{\b} + s$ vertices and at least
$2s$ edges.  An Euler characteristic argument then shows that $H$ has
a disk face $E$.  This disk $E$ contains a 1-cycle of $G_{\b}$ (see
[CGLS, p.\ 279] for definition), and hence a Scharlemann cycle [CGLS,
Lemma 2.6.2].  But this contradicts the fact that $E$ is a face of
$H$, because by definition $H$ would have a vertex in the disk bounded
by this Scharlemann cycle.  \endproof

Let $P, Q$ be two edge endpoints on the boundary of a vertex $u$ in
$G_{\a}$.  Let $P_0 = P, P_1, \ldots, P_{k-1}, P_k = Q$ be the edge
endpoints encountered when traveling along $\bdd u$ in the direction
induced by the orientation of $u$.  Then the {\it distance\/} from $P$
to $Q$ (at the vertex $u$) is defined as $\rho_u (P,Q) = k$.  Notice
that since the valency of $u$ is $\Delta n_{\b}$, we have $\rho_u(Q,
P) =\Delta n_{\b} - \rho_u (P, Q)$.  If $e_1, e_2$ is a pair of edges,
each having a single endpoint $P_i$ on the vertex $u$ in $G_{\a}$,
then define $\rho_u(e_1, e_2) = \rho_u (P_1, P_2)$.

A pair of edges $e_1, e_2$ connecting two vertices $u, v$ in $G_{\a}$
is an {\it equidistant pair\/} if $\rho_u(e_1, e_2) =\rho_v(e_2,
e_1)$.  In particular, one can check that if $e_1, e_2$ are a pair of
parallel edges connecting a pair of parallel vertices in $G_{\a}$,
then $e_1, e_2$ is an equidistant pair in $G_{\a}$.  

\proclaim{Lemma 2.4} {\rm (The Equidistance Lemma.)}  Let $e_1, e_2$
be a pair of edges connecting the same vertices on $G_1$ and the same
vertices on $G_2$.  Then $e_1, e_2$ is an equidistant pair in $G_1$ if
and only if it is an equidistant pair in $G_2$.  \endproclaim

\proof
See [GW1, Lemma 2.8].
\endproof

Given two slopes $r_1, r_2$ on the torus $T_0$, let $l$ be a curve
intersecting $r_1$ at a single point.  Choosing $l$ and the
orientations of the curves properly, we may assume that homologically
$r_2 = q r_1 + \Delta l$, where $1\leq q < \Delta/2$.  The number $q$
is called the {\it jumping number\/} of $r_1, r_2$.  Note that if
$\Delta = 4$ then $q = 1$, and if $\Delta = 5$ then $q=1$ or $2$.

\proclaim{Lemma 2.5} (1) If the jumping number $q=1$, in particular if
$\Delta = 4$, then a pair of $j$-edges at a vertex $u_i$ in $G_{\a}$
are adjacent among all the $j$-edges if and only if on $G_{\b}$ they
are also adjacent at $v_j$ among all $i$-edges.  

(2) If $q=2$, then a pair of $j$-edges at a vertex $u_i$ in $G_{\a}$
are adjacent among all $j$-edges if and only if on $G_{\a}$ they are
not adjacent among all the $i$-edges at $v_j$.  \endproclaim

\proof
This is essentially [GW1, Lemma 2.10].  It was shown that if $P_1,
..., P_{\Delta}$ are the consecutive $j$-edge endpoints at $u_i$, then
on $\bdd v_j$ they appear in the order $P_q, P_{2q}, ..., P_{\Delta
q}$, hence the result follows.
\endproof

A graph $G$ on an annulus $A$ is {\it special\/} if every vertex has
at least two nonparallel boundary edges.  Note that $G$ is special if
and only if the corresponding reduced graph $\hat G$ is special.

\proclaim{Lemma 2.6} 
(1) If $G$ is special then every vertex has exactly two boundary edges
in $\hat G$, going to distinct boundary components of $A$.

(2) If $G_{\a}$ has $2n_{\b}$ parallel boundary edges, then $G_{\b}$
is special.  $G_{\a}$ cannot have more than $2n_{\b}$ parallel
boundary edges.

(3) If some edge $E$ of $\hat G_{\a}$ represents $n_{\b}$ negative
edges, and if $G_{\a}$ has some positive edges, then $G_{\a}$ has at
most $n_{\b}$ parallel boundary edges, and each vertex of $\hat
G_{\b}$ has at most one boundary edge.  \endproclaim

\proof
(1) Otherwise there would be a pair of edges of $\hat G$ at some 
vertex $v$ 
going to the same boundary component of $A$.  By looking at an 
outermost such pair one can see that some vertex $u$ of $\hat G$ has a 
single boundary edge in $\hat G$, contradicting the definition of 
a special graph.

(2) If $G_{\a}$ has $2n_{\b}$ parallel boundary edges, then for any
label $i$ it has two parallel $i$-edges.  Since no two edges are
parallel on both graphs, these two edges are non-parallel on $G_{\b}$,
hence $G_{\b}$ is special.  If $G_{\a}$ has more than $2n_{\b}$
parallel boundary edges, then there is a label $i$ such that $G_{\a}$
has three parallel boundary $i$-edges.  Since by (1) the vertex $v_i$
in $\hat G_{\b}$ has only two boundary edges, two of these edges would
be parallel on both graphs, contradicting Lemma 2.1(2).

(3)  Since $G_{\a}$ has some positive edges, the vertices of $G_{\b}$
cannot all be parallel, so there are at least two $E$-orbits, which 
form parallel essential cycles on $G_{\b}$.  Hence all boundary edges 
at a vertex of $G_{\b}$ must be parallel to each other.  If $G_{\a}$ has
more than $n_{\b}$ parallel boundary edges then two of them would be
parallel on both graphs, contradicting Lemma 2.1(2).
\endproof

\proclaim{Lemma 2.7} Suppose all vertices of $\hat G_{\a}$ are
boundary vertices, and suppose there are two boundary edges $E_1, E_2$
of $\hat G_{\a}$ incident to the same vertex $v$ and going to the same
boundary component of $A_{\a}$.  Then $\hat G_{\a}$ has a vertex $v'$
of valency at most 3 which is incident to a single boundary edge,
and $G_{\b}$ is special.  \endproclaim

\proof Let $D$ be the disk on $A_{\a}$ cut off by $E_1\cup E_2$.
Since $E_1, E_2$ are nonparallel, $D$ contains a vertex $v_1 \neq v$,
hence by adding an edge if necessary we may assume that there is an
edge incident to $v$ other than $E_1, E_2$.  Let $\tilde D$ be the
double of $D$ along $E_1\cup E_2$, and let $\tilde G$ be the double of
$\hat G_{\a} \cap D$.  Then each vertex of $\tilde G$ has a boundary
edge.  By [CGLS, Lemma 2.6.5] $\tilde G$ has a vertex $v'$ of valency
at most 3 and incident to at most one boundary edge.  Since $v$
has valency at least $4$ in $\tilde G$, $v' \neq v$.  Hence $\val (v',
\hat G_{\a}) = \val (v', \tilde G) \leq 3$.  By Lemma 2.2(3) each
interior edge of $\tilde G_{\a}$ represents at most $n_{\b}$ edges.
Since $\Delta \geq 4$, this implies that the unique boundary edge at
$v'$ represents at least $2n_{\b}$ edges.  By Lemma 2.6(2) in this
case $G_{\b}$ is special.  \endproof

\head  \S 3.  Reduced graphs on annuli      \endhead

By a {\it reduced graph\/} on a surface we mean one with no trivial
loops or parallel edges; in other words, no faces of the graph are
monogons or bigons.

\definition{Definition 3.1}  Let $G$ be a reduced graph on an annulus
$A$.  Then $G$ is said to be {\it triangular\/} if 

(i)  every vertex has at most one boundary edge; 

(ii)  every interior vertex has valency 6;

(iii)  every boundary vertex has valency 5;

(iv)  every face of $G$ is a disk with three edges.
\enddefinition

We remark that the only properties of a triangular graph that we will
use are (i), (iii), and the fact that the graph has at least one
boundary vertex (which follows from (iv)).

\proclaim{Proposition 3.1}  Let $G$ be a reduced graph in an annulus $A$.
Then either 

(1)  $G$ contains an interior vertex of valency at most $5$; or

(2)  $G$ contains a boundary vertex of valency at most $4$ with
exactly one boundary edge; or 

(3)  $G$ is triangular; or

(4)  $G$ is special.
\endproclaim

\proof Let $G_1$ be a graph obtained from $G$ by adding extra edges so
as to make each face of $G_1$ a disk with three edges.  In particular,
in $G_1$, each boundary face has three edges, and if some vertex $v$
has two boundary edges $e_1, e_2$, then $v$ has an edge on each side
of $e_1 \cup e_2$, so it has valency at least 4.

Let $G_2$ be the union of $G_1$ and $\bdd A$, with the obvious graph
structure.  Thus the points of $G_1\cap \bdd A$ are now considered
vertices, and the segments of $\bdd A$ cut by these vertices are
considered edges of $G_2$.  Note that $\val (v, G_2) = 3$ for all
vertices $v$ on $\bdd A$, and each boundary face now has four edges.
Let $G_3$ be obtained from $G_2$ by adding a diagonal edge in each
boundary face of $G_2$, all sloping in the same direction; in other
words, no two edges added have a common vertex on $\bdd A$.  We have
$\val (v, G_3) = 4$ if $v \in \bdd A$.  One can see that if we remove
all edges and vertices on $\bdd A$ then we get a graph that is
obtained from $G_1$ by adding an extra copy of each boundary edge.
Hence if $\val (v, G_1) = p$ and $v$ has $q$ boundary edges in $G_1$,
then $\val(v, G_3) = p+q$.  In particular, if $v$ has two boundary
edges in $G_1$, then its valency in $G_3$ is at least $4+2 = 6$.  Each
face of $G_3$ is now a triangle.

The double of $A$ along $\bdd A$ is a torus $T$, and the corresponding
double of $G_3$ is a reduced graph $\tilde G_3$ on $T$ with triangular
faces.  By an Euler characteristic argument, one can show that the
number of edges in $\tilde G_3$ is three times the number of vertices
of $\tilde G_3$.  Thus either (i) some vertex $v$ of $\tilde G_3$ has
valency at most 5, or (ii) all vertices of $\tilde G_3$ have valency
6.  All vertices on $\bdd A$ have valency 4 in $G_3$, hence valency 6
in $\tilde G_3$, and we have shown that if $v$ has two boundary edges
in $G_1$ then it has valency at least $6$ in $G_3$; therefore (i)
implies that either $v$ is an interior vertex of $G$ with valency at
most 5, or it is a boundary vertex of $G$ with valency at most $5 - q
\leq 4$ and incident to at most one boundary edge, so the graph is of
type (1) or (2) in the proposition.  Hence we may assume that all
vertices of $G_3$ in the interior of $A$ have valency 6.

If no vertex of $G_1$ has two boundary edges then each boundary vertex
of $G_1$ has valency $6 - q = 6 - 1 = 5$.  Since each interior vertex
of $G_1$ has valency 6, it follows that $G_1$ is triangular.  Since
$G$ is a subgraph of $G_1$ with the same vertices, either $G = G_1$
and hence $G$ is of type (3), or $G$ has a vertex $v$ with $\val (v,
G) < \val (v, G_1)$, in which case $G$ is of type (1) or (2).

Now assume some vertex $v$ of $G_1$ has two boundary edges $e_1, e_2$
going to different boundary components.  Then the valency of $v$ in
$G_1$ is at most $6 - 2 = 4$.  Since each face of $G_1$ has three
edges, there is exactly one interior edge $e'$ on each side of
$e_1\cup e_2$.  Let $v'$ be the other endpoint of $e'$.  Since each
face has three edges, $v'$ must also have two boundary edges going to
different boundary components of $A$.  Repeating this process, we see
that $G_1$ is a special graph such that each vertex has valency 4.
Since $G$ is a subgraph of $G_1$, either it is special, hence of type
(4), or it has a vertex of valency at most 3 and incident to at most
one boundary edge, in which case it is of type (1) or (2).

Finally, assume $G_1$ has a vertex $v$ which has two boundary edges
going to the same boundary component.  Then they cut off a disk $D$
from the annulus, which we may assume to be outermost.  However,
arguing as in the previous paragraph, we see that the vertex on the
other end of an edge $e'$ in $D$ incident to $v$ must have two
boundary edges, which is a contradiction.  Therefore this case does
not happen.  \endproof

\head   \S 4.   Special graphs
\endhead

Recall that a graph $G$ on an annulus $A$ is special if every vertex
has two nonparallel boundary edges.  By Lemma 2.6(1) this implies that
every vertex of $G$ has exactly two boundary edges in $\hat G$, going
to different boundary components of $A$.

To simplify notation, denote $n_{\b}$ by $n$.

\proclaim{Lemma 4.1}  
If $G_{\a}$ is a special graph, then $G_{\b}$ is also special.
\endproclaim

\proof First notice that since each vertex $u_i$ of $G_{\a}$ is
incident to at most two families of interior edges and each such
family contains at most $n$ edges (Lemma 2.2(3)), there are at most
two interior $j$-edges at $u_i$ for any $j$.  Hence there are at least
$\Delta -2$ boundary $j$-edges at $u_i$.  Since this is true for all
$i, j$, we see that each vertex $v_j$ of $G_{\b}$ has at least $2
n_{\a}$ ($3n_{\a}$ if $\Delta = 5$) boundary edges.  Since each
parallel family contains at most $2n_{\a}$ edges (Lemma 2.6(2)), the
lemma follows immediately when $\Delta = 5$.

Now assume $\Delta = 4$, and assume $G_{\b}$ is not special.  Then it
has a vertex $v_i$ such that all boundary edges are parallel.  By
Lemma 2.6(2) and the above, $v_i$ has exactly $2n_{\a}$ boundary
edges, all parallel to each other.  In particular, there are only two
boundary $1$-edges $e_1, e'_1$ at $v_i$.  Dually this means that $e_1,
e'_1$ are the only boundary $i$-edges at $u_1$.  Since they are
parallel on $G_{\b}$, they cannot be parallel on $G_{\a}$, so they
belong to different families of boundary edges.  Since these two edges
are adjacent among all $1$-edges at $v_i$, by Lemma 2.5(1) they must
also be adjacent among all $i$-edges at $u_1$.  This implies that the
two interior $i$-edges at $u_1$ are on the same side of $e_1\cup
e'_1$, so they belong to the same edge $E$ in $\hat G_{\a}$ because
there is only one interior edge of $\hat G_{\a}$ on each side of
$e_1\cup e'_1$.  Since by Lemma 2.2(3) $E$ contains at most $n$
edges, this is impossible.  \endproof

In the remainder of this section we will assume that both $G_1$ and
$G_2$ are special.

The sign of a vertex $u$ in $G_{\a}$ induces an orientation on $\bdd
u$, called its {\it preferred orientation}.  Thus the preferred
orientations of the $\bdd u$'s are all in the same direction on $T_0$.
Let $e_1, e_2$ be a pair of adjacent boundary edges at some vertex $u$
of $G_{\a}$.  When traveling on $\bdd u$ along the preferred
orientation, the labels at the endpoints of $e_1, e_2$ appear as $i,
i+1$ for some $i$ ($i+1=1$ if $i = n$).  They cut off a band $B$ on
the surface $F_{\a}$, called an {\it $i$-band} at $u$ (of $G_{\a}$).
Note that the label $i$ is determined by the pair $e_1, e_2$ even if
$n = 2$.  The edge labeled $i$ at $u$ is called the {\it initial
edge\/} of $B$, the other the {\it terminal edge.}  Two $i$-bands of
$G_{\a}$ are of different {\it types\/} if their initial edges are
nonparallel on $G_{\b}$; otherwise they are of the same type.

If $e_1, ..., e_k$ are all the edges of a parallel family $E$ at a
vertex $u$, appearing in this order when traveling along the preferred
orientation of $\bdd u$, then $e_k$ is called the {\it ending edge\/}
of $E$, and the label of $e_k$ at $u$ is called the {\it ending
label\/} of $E$.  Note that if a boundary $i$-edge $e$ is not an
ending edge, then it is the initial edge of an $i$-band.

\proclaim{Lemma 4.2}  There is a label $i$ such that all $i$-bands of
$G_{\a}$ are of the same type.
\endproclaim

\proof Assuming otherwise, then there are two $i$-bands $B_i^1, B_i^2$
of different types for each $i$.  Since the graph $G_{\b}$ is special,
there are only two families of parallel boundary edges for each vertex
$v_i$ in $G_{\b}$, so each family contains the initial edge of some
$B_i^j$.  Therefore, the terminal edge of each $B_i^j$ is parallel to
the initial edge of some $B_{i+1}^k$, so there is a band $D_i^j$ on
$F_{\b}$ connecting these two edges.  Note that $D_i^j$ degenerates to
a single edge if these two edges coincide.  

Consider the 2-complex $Q = \cup (B_i^j \cup D_i^j)$.  Then $Q \cap
T_0 = \cup (e_i^j \cup d_i^j)$ is a graph $G$ on $T_0$, where $e_i^j =
B_i^j \cap T_0$ and $d_i^j = D_i^j \cap T_0$.  We have $Q \cong G
\times I$.  Shrinking each $d_i^j$ to a point, and orienting $e_i^j$
so that its endpoint is on $d_i^j$, we get an oriented graph $G'$ in
which each vertex $d_i^j$ is the tail of some edge $e_{i+1}^k$.  Hence
$G'$ contains an embedded oriented cycle.  The corresponding cycle $C$
in $G$ is then an embedded loop in $T_0$.  Let $\gamma$ be a parallel
copy of some boundary component of $F_{\b}$ on $T_0$, intersecting
some $e_i^j$ in $C$ transversely at a single point.  The definition of
$B_i^j$ and the orientation of $e_i^j$ implies that $C$ intersects
$\gamma$ always in the same direction; hence $C$ is an essential
curve.  Thus $A = C\times I \subset Q$ is an annulus properly embedded
in $M$ intersecting $T_0$ in the essential curve $C$, which
contradicts the assumption that $M$ is $\bdd$-irreducible and
anannular.  \endproof

\proclaim{Lemma 4.3}  Each family of boundary edges in $G_{\a}$
contains at least $n$ edges.  \endproclaim

\proof Let $E_1, ..., E_4$ be the four edges of $\hat G_{\a}$ at
$u_1$, with $E_1, E_3$ the boundary edges.  If $|E_1| < n$ then there
is a label $i$ which does not appear at the endpoints of edges in
$E_1$.  If $\Delta = 5$ then we would have $|E_3| = 5n - |E_1| - |E_2|
- |E_4| \geq 2n + 1$, contradicting Lemma 2.6(2).  Hence $\Delta = 4$.
Since $|E_3| \leq 2n$, $E_3$ contains at most two $i$-edges, so each
of $E_2, E_4$ contains one $i$-edge.
Let $e_1, e_2$ be the $i$-edges of $E_2, E_4$ at $u_1$, and let $e_3,
e_4$ be the $i$-edges of $E_3$.  Since $e_3, e_4$ are adjacent
$i$-edges at $u_1$, by Lemma 2.5(1) they are adjacent $1$-edges at
$v_i$.  On $G_{\b}$ the two edges $e_3, e_4$ belong to different
families of boundary edges at $v_i$, because they cannot be parallel
on both graphs.  Therefore the two edges $e_1, e_2$ belong to the same
family of interior edges.  Since they both have label $1$ at $v_i$,
this would imply that the interior family containing them has
at least $n_{\a} + 1$ edges, contradicting Lemma 2.2(3).  \endproof

\proclaim{Lemma 4.4}  The jumping number $q = 1$.
\endproclaim

\proof This is automatically true if $\Delta = 4$.  Hence assume
$\Delta = 5$.  First assume that there is a vertex $u_j$ of $G_{\a}$
which has two interior $i$-edges $e_1, e_2$ for some $i$.  Since each
interior family contains at most $n$ edges, $e_1, e_2$ are nonparallel
on $G_{\a}$.  By Lemma 4.3 each boundary family contains an $i$-edge,
hence $e_1, e_2$ are non adjacent among the $i$-edges at $u_j$.
Dually on $G_{\b}$ these are $j$-edges at the vertex $v_i$.  For the
same reason, they are non adjacent among all $j$-edges at $v_i$.
Therefore by Lemma 2.5 the jumping number $q =1$.

Now assume that $u_j$ has at most one interior $i$-edge for all $i$.
Then it has at most $n$ interior edge endpoints.  On the other hand,
since each boundary family contains at most $2n$ edges and the valency
of $u_j$ is $\Delta n = 5n$, we see that it cannot have less than $n$
interior edge endpoints; therefore it has exactly $n$ interior edge
endpoints, and each boundary family contains exactly $2n$ edges.  If
$u_j$ has two interior families, so each family contains less than $n$
edges, then the two boundary families have different ending labels.
In this case for each label $i$ there are three $i$-bands, which
cannot all be of the same type because each boundary family of $v_i$
has at most two $j$-edges.  This contradicts Lemma 4.2.  Therefore
$u_j$ has only one family of interior edges, which contains $n$ edges.
For the same reason, each vertex of $G_{\b}$ has only one family of
interior edges, containing $n_{\a}$ edges.  By the parity rule one of
these families is negative, and by Lemma 2.1(3) they form cycles on
the other graph, so each vertex of that graph would then have two
families of interior edges, contradicting the above conclusion.
\endproof

\proclaim{Lemma 4.5}  Suppose all $i$-bands at a vertex $u_j$ of
$G_{\a}$ are of the same type.  Then 

(1) there are $n$ parallel interior edges at $u_j$, and

(2) each family of $n$ parallel interior edges at $u_j$ has $i$
as its ending label.  \endproclaim

\proof Let $E_1, ..., E_4$ be the edges at $u_j$ of $\hat G_1$,
appearing in this order around $\bdd u_j$ along its preferred
orientation, with $E_1, E_3$ the boundary edges.  Let $e_1, ..., e_4$
be four $i$-edges at $u_j$, appearing successively along the preferred
orientation of $\bdd u_j$.

First assume that all $e_i$ are boundary edges.  Then we may assume
that $e_1, e_2 \in E_1$, and $e_3, e_4 \in E_3$.  Thus $e_1, e_3$ are
not ending edges, so they are initial edges of some $i$-bands $B_1,
B_3$.  Since the jumping number $q=1$ (Lemma 4.4), and since $e_1,
e_3$ are non adjacent among $i$-edges at $u_j$, by Lemma 2.5 they are
non adjacent among $1$-edges at $v_i$ in $G_{\b}$, hence they are non
parallel boundary edges on $G_{\b}$.  Therefore $B_1, B_3$ are of
different type.

Now assume that $E_2$ contains an $i$-edge $e_2$, say.  Since each of
$E_1, E_3$ contains at least $n$ edges, we must have $e_1 \in E_1$ and
$e_3 \in E_3$.  Assume that either $e_2$ is not the ending edge of
$E_2$ or $|E_2| < n$.  Then $e_1$ is not an ending edge of $E_1$, and
there is an $i$-band $B_1$ with $e_1$ as the initial edge.  If $e_3$
is not an ending edge either, then there is an $i$-band $B_3$ with
$e_3$ as initial edge.  For the same reason as above, $B_1, B_3$ are
of different type, and we are done.  So assume that $e_3$ is the
ending edge of $E_3$.  Now we must have $|E_4| < n$ as otherwise $E_4$
would have $n$ edges and have the $i$-edge $e_4$ as its ending edge,
contradicting the assumption.  Hence $e_4$ is in $E_1$, and so there
is an $i$-band with $e_4$ as an initial edge.  Since $e_1, e_4$ are
parallel on $G_{\a}$, they are nonparallel on $G_{\b}$, so again $B_1,
B_4$ are of different type.  This proves (1).  To prove (2), notice
that if $|E_2| = n$ but $e_2$ is not the ending edge, then $e_1, e_3$
are not ending edges of $E_1, E_3$, so from the above the two $i$-bands
$B_1, B_3$ are of different type.  \endproof

\proclaim{Lemma 4.6} Suppose $n > 2$.  Then each positive edge of
$\hat G_{\a}$ represents at most $n/2$ edges.  \endproclaim

\proof When $n > 2$, the special graph $\hat G_{\b}$ has at most
one edge connecting any two vertices.  If $G_{\a}$ has $n/2+1$
parallel positive edges, then it has a Scharlemann cycle $e_1\cup e_2$
with label pair $(1,2)$, say.  So the two edges $e_1, e_2$ would be
parallel on $G_{\b}$, contradicting Lemma 2.1(4).  \endproof

\proclaim{Proposition 4.7} If $G_{\a}$ is special then up to
relabeling we have $n_1=1$, $n_2 = 2$, and $G_1$ has exactly two
interior edges.  \endproclaim

\proof First assume $n_{\a}\geq 2$ for $\a = 1,2$.  By Lemma 4.2 for
each graph $G_{\a}$ there is a label $i$ such that all $i$-bands of
$G_{\a}$ are of the same type.  Let $u_j$ be a vertex of $G_{\a}$.  By
Lemma 4.5(1), it has a set of $n$ parallel interior edges $E$ with $i$
as its ending label at $u_j$.  Let $u_k$ be the vertex on the other
endpoint of $E$, then by Lemma 4.5(2), $E$ also has ending label $i$
at $u_k$.  If $E$ is negative, then the ending edge $e$ of $E$ at
$u_j$ is the same as that at $u_k$, so $e$ would have the same label
$i$ on its two endpoints, and hence is a loop on $G_{\b}$.  Since $n
\geq 2$ and $G_{\b}$ is special, this is absurd.  If $E$ is positive,
then the two ending edges would give rise to two negative edges at
$v_i$ in $G_{\b}$, which must be nonparallel because they cannot be
parallel on both graphs.  Thus both families of interior edges at
$v_i$ are negative.  Replacing $u_j$ by $v_i$ in the above argument,
we get a contradiction because now $E$ must be negative.  Therefore up
to relabeling we must have $n_1 = 1$.

By Lemma 4.5(1) the only vertex $u_1$ of $G_1$ has $n_2$ parallel
interior edges, which by the parity rule must be negative edges on
$G_2$, hence $n_2 \geq 2$.  If $n_2 > 2$, then by Lemma 4.6 $G_1$ has
at most $n_2/2$ interior edges, which is a contradiction.  Hence the
result follows.  \endproof

\head  \S 5.  The generic case
\endhead

In this section except for Lemma 5.1, we assume $n_{\a}, n_{\b} > 2$.
By Proposition 4.7, $G_{\a}$ and $G_{\b}$ are not special.  Again
denote $n_{\b}$ by $n$.

\proclaim{Lemma 5.1} Suppose $n_{\a} \geq 2$, $n > 2$, and suppose
$\hat G_{\a}$ has a negative edge $E$ with $|E| = n$, and a positive
edge $E'$ with $|E'| = n/2+1$.  Let $(1,2)$ be the label pair of the
Scharlemann cycle in $E'$.  Then

(1) $G_{\b}$ has at most $n/2$ boundary vertices; 

(2) when $n=4$, the two vertices $v_3, v_4$ of $G_{\b}$ cannot
both be boundary vertices; and

(3) when $n=4$, $G_{\a}$ cannot have both a $(1,4)$-edge and a
$(2,3)$-edge.  \endproclaim

\proof Let $k$ be the number of $E$-orbits.  Since $E'$ contains more
than $n/2$ edges, hence contains a Scharlemann cycle, the annulus
$A_{\b}$ is separating, so $G_{\b}$ has the same number of positive
and negative vertices.  Each $E$-orbit contains the same number
($n/k$) of vertices, all of the same sign, so the number of orbits
containing positive vertices is the same as the number of those
containing negative ones, and hence $k$ must be even.  Recall that
each $E$-orbit forms an essential cycle on $G_{\b}$, so only the
vertices on the two cycles adjacent to the two boundary components of
$A_{\b}$ could be boundary vertices.  Hence the number of boundary
vertices is at most $2(n/k)$, and since $k$ is even, (1) follows 
unless $k = 2$.

Assume $k=2$.  Let $C_1, C_2$ be the two cycles of $E$-orbits on
$G_{\b}$, and let $e_1, e_2$ be the edges of the Scharlemann cycle in
$E'$.  By Lemma 2.1(4) $e_1\cup e_2$ is an essential cycle on
$G_{\b}$.  The two vertices $v_1, v_2$ of $e_1\cup e_2$ are on
different $C_1, C_2$ because they are antiparallel, so the cycle
$e_1\cup e_2$ lies between $C_1$ and $C_2$, separating the vertex
$v_3$ on the first orbit from the vertex $v_n$ on the second.  On the
other hand, since $E'$ contains more than two edges, there is an edge
adjacent to the Scharlemann cycle which is a $(3, n)$-edge, so on
$G_{\b}$ there would be an edge connecting $v_3$ to $v_n$.  This is a
contradiction, showing that $k=2$ is impossible.  In particular, this
proves (1).

Now assume $n=4$.  Since we have shown that $k$ is even and $k\neq 2$,
we must have $k=4$.  In this case each vertex $v_i$ of $G_{\b}$ has an
essential loop $C_i$ coming from the $n$ parallel negative edges in
$G_{\a}$.  These loops and their vertices form essential circles on
$A_{\b}$ which are parallel to each other.  As above, there is an edge
in $E'$ which connects $v_3$ to $v_4$.  Hence the circles $C_3$ and
$C_4$ are adjacent to each other, so $v_3, v_4$ cannot both be
boundary vertices.  This proves (2).  Since the edges in the
Scharlemann cycle connect $v_1$ to $v_2$, $C_1$ is adjacent to $C_2$.
Thus either $C_3$ separates $v_4$ from $v_1, v_2$, so there is no edge
connecting $v_4$ to $v_1$, or $C_4$ separates $v_3$ from $v_1,v_2$, so
there is no edge connecting $v_3$ to $v_2$.  This proves (3).
\endproof

\proclaim{Lemma 5.2} Suppose $E_1, ..., E_5$ are the edges of $\hat
G_{\a}$ at a vertex $u$ of valency 5.  If $E_1, E_2, E_3$ are
positive, and $E_4$ is an interior edge, then $|E_5| > n$; in
particular, $E_5$ is a boundary edge.  \endproclaim

\proof Assume $|E_5| \leq n$.  By Lemma 2.2(2) we have $|E_i| \leq n/2
+ 1$ for $i \leq 3$, and by Lemma 2.2(3) $|E_4| \leq n$.  Since $n>2$,
we must have $\Delta = 4$, and
$$ 4n = \Delta n = |E_1| + \cdots + |E_5| \leq 3 (\frac {n}{2} + 1) +
2n = \frac 72 n + 3,$$ which implies that $n \leq 6$.  Moreover, $n$
must be even, otherwise by Lemmas 2.2(2) and 2.1(4) we would have
$|E_i| \leq n/2$ for $i=1,2,3$, hence $ 4n \leq 3 (n/2)+ 2n $, which
is absurd.

If $n = 6$, then all the above inequalities are equalities.  In
particular, $\Delta = 4$, $|E_4|=|E_5|= 6$, and $|E_i| = 4$ for $i =
1, 2,3$, so each of $E_1, E_2, E_3$ contains a Scharlemann cycle, and
by Lemma 2.2(1) they all have the same label pair, say $(1,2)$.
But since $|E_4|= | E_5| = n$, these labels also appear in $E_4$ and
$E_5$.  Thus the label $1$ appears 5 times, contradicting the fact
that $\Delta = 4$.

Now assume $n = 4$.  If each of $E_1, E_2, E_3$ contains a Scharlemann
cycle with label pair $(1,2)$, say, (in particular, if $|E_i| = 3$
for $i=1,2,3$), then again the labels $\{1,2\}$ appear three times
among the endpoints of $E_1\cup E_2 \cup E_3$ at $u$.  Also, since
$E_4\cup E_5$ has at least $16 - 3\times 3 = 7$ edge endpoints at $u$,
one of the labels $\{1,2\}$ appears at least twice among the endpoints
of $E_4\cup E_5$ at $u$, so it appears 5 times at $u$, contradicting
the fact that $\Delta = 4$.  Hence we may assume that $|E_1| = |E_2| =
3$, $|E_3| = 2$, $E_3$ contains no Scharlemann cycle, and $|E_4|=
|E_5| = 4$.  Since the two edges of $E_3$ have labels $3, 4$ at $u$,
they must have label sets $\{1,4\}$ and $\{2,3\}$.  Since $|E_4| = 4$,
the edges in $E_4$ are negative.  This contradicts Lemma 5.1(3),
completing the proof of the lemma.  \endproof

\proclaim{Lemma 5.3}  $\hat G_{\a}$ has no interior vertex of valency
at most $5$.
\endproclaim

\proof Let $E_1, \ldots, E_5$ be the edges of $\hat G_{\a}$ incident
to $u$.  Since all these edges are interior edges, by Lemma 5.2 they
can have at most two positive edges, say $E_1, E_2$.  By Lemma
2.3(2), $u$ has at most $2n$ negative edges in $G_{\a}$, hence $E_1
\cup E_2$ represents at least $2n$ positive edges.  By Lemma 2.2(2)
we have $2n \leq 2(n/2+1)$, which contradicts the assumption that $n
\geq 3$.  \endproof

\proclaim{Lemma 5.4} $\hat G_{\a}$ cannot have a boundary vertex $u$
of valency at most 4 with a single boundary edge.  \endproclaim

\proof Let $E_0$ be the boundary edge, and $E_1, E_2, E_3$ the
interior edges of $\hat G_{\a}$ at $u$.  By Lemma 2.6(2) and
Proposition 4.7 we have $|E_0| < 2n$.  By Lemma 2.3(2), $u$ can have
at most $2n$ negative edges in $G_{\a}$, so one of the interior edges,
say $E_1$, must be positive, and by Lemma 2.2(2) $|E_1| \leq n/2+1 <
n$.  Since each of $E_2, E_3$ represents at most $n$ edges, we have
$|E_0| > n$.

We claim that {\it either $G_{\a}$ or $G_{\b}$ contains a Scharlemann
cycle.}  If $u$ has more than $n$ negative edges, then by Lemma
2.3(1)  $G_{\b}$ contains a Scharlemann cycle.  So assume $u$ has at
most $n$ negative edges.  Since $u$ has less than $2n$ boundary edges,
it must have more than $n$ positive edges.  If at most two of the
$E_1, E_2, E_3$ are positive, then one of them represents more than
$n/2$ positive edges; if all the three interior edges at $u$ are
positive, then since they represent more than $2n$ edges, again one
of them represents more than $n/2$ edges.  In either case these
parallel edges contain a Scharlemann cycle.  This completes the proof
of the claim.

Now $|E_0| > n$ implies that some vertex of $G_{\b}$ has two
nonparallel boundary edges.  In particular, $\hat G_{\b}$ cannot be
triangular.  It follows from Lemma 5.3 and Proposition 3.1 that $\hat
G_{\b}$ must also have a boundary vertex $v$ of valency at most 4 with
a single boundary edge.  Since one of $G_{\a}$ and $ G_{\b}$ has a
Scharlemann cycle, by considering $v$ instead of $u$ if necessary, we
may assume without loss of generality that $G_{\a}$ contains a
Scharlemann cycle with label pair $(1, 2)$.

We claim that $|E_0| \leq n+2$.  Otherwise each vertex of $G_{\b}$ has
a boundary edge, and some vertex $v_i$ other than $v_1, v_2$ has two
such edges $e_1, e_2$.  Since the edges of the Scharlemann cycle form
an essential subgraph of $G_{\b}$ (Lemma 2.1(4)), separating the two
boundary components of $A_{\b}$, the edges $e_1, e_2$ must go to the
same boundary component.  Applying Lemma 2.7, we see that $G_{\a}$ is
special, a contradiction.

Since $|E_0| > n$ and $u$ has some positive edges, by Lemma 2.6(3) the
graph $G_{\a}$ cannot have $n$ parallel negative edges.  Thus if $k$
of the $E_1, E_2, E_3$ are positive, then
$$ 4n \leq (n+2) + k (\frac n2 +1) + (3-k)(n-1) = (4-\frac k2)n + (2k
-1)$$ which implies that $n < 4$.  But since $G_{\a}$ contains a
Scharlemann cycle, $n$ is even.  This contradicts the assumption that
$n > 2$.  \endproof

\proclaim{Lemma 5.5} If both $\hat G_1, \hat G_2$ are triangular, then
each boundary vertex has exactly two positive and two negative edges
in $\hat G_{\a}$.  \endproclaim

\proof Let $E_0, ..., E_4$ be the edges of $\hat G_{\a}$ at a boundary
vertex $v$, with $E_0$ the boundary edge.  Since $\hat G_{\b}$ is also
triangular, $E_0$ represents at most $n$ edges.  Therefore by Lemma
5.2 at most two of the $E_i$ are positive.  On the other hand, by
Lemma 2.3(2) $v$ has at most $2n$ negative edges, hence at least $n$
positive edges.  Since each $E_i$ represents at most $n/2 + 1 < n$
positive edges, $v$ must have two positive edges in $\hat G_{\a}$.
\endproof

\proclaim{Lemma 5.6} Suppose both $\hat G_1, \hat G_2$ are triangular.
Then all vertices of $G_1, G_2$ are boundary vertices.  \endproclaim

\proof Let $E_0, ..., E_4$ be the edges of $\hat G_{\a}$ at a boundary
vertex $v$, with $E_0$ the boundary edge.  By Lemma 5.5 we may
assume that $E_1, E_2$ are negative edges, and $E_3,E_4$ are positive
edges.

Suppose $G_{\b}$ has some interior vertices.  Then $|E_0| < n$.  Since
$v$ has at most $2n$ negative edges, it has more than $n$ positive
edges, so $E_3\cup E_4$ contains a Scharlemann cycle with label pair
$(1,2)$, say, and $n$ is even.  Also, one of $E_1, E_2$ must
represent $n$ parallel edges, for otherwise $E_0\cup E_1 \cup E_2$
would contain at most $3n-3$ edges, so one of $E_3, E_4$ would
contain at least $n/2+2$ edges, contradicting Lemma 2.2(2).  Now we
can apply Lemma 5.1(1) and conclude that $G_{\b}$ has at most $n/2$
boundary vertices.  Thus $|E_0| \leq n/2$.  We have the inequality
$$ 4n \leq |E_0| + ... + |E_5| \leq \frac n2 + 2n + 2(\frac n2 + 1) \leq
\frac 72 n + 2.$$ Since $n$ is even, this implies $n=4$, $|E_0| = 2$,
$|E_1|=E_2|=4$, and $|E_3| = |E_4| = 3$.  Now each of $E_3, E_4$
contains a Scharlemann cycle on label pair $(1,2)$, so these labels
appear 4 times among the interior edge endpoints at $v$.  Thus the
labels of $E_0$ must be $3,4$.  This contradicts Lemma 5.1(2).
\endproof

\proclaim{Lemma 5.7}  $\hat G_1, \hat G_2$ cannot both be triangular.
\endproclaim

\proof Assume $\hat G_1, \hat G_2$ are triangular.  Then by Lemma 5.6
all vertices of $G_1, G_2$ are boundary vertices, and by Lemma 5.5
each vertex $v$ of $G_{\a}$ has exactly two positive edges and two
negative edges in $\hat G_{\a}$.  Since a positive edge in $G_{\a}$ is
a negative edge in $G_{\b}$, it follows that either (i) some vertex
$v$ of one of the graphs, say $G_1$, has more positive edge endpoints
than negative ones, or (ii) all vertices of $G_1$ and $G_2$ have the
same number of positive and negative edge endpoints.

In case (i), (writing $n = n_2$), $v$ has at most $2(n/2+1) = n + 2$
positive edges, at most $n + 1$ negative edges, and at most $n$
boundary edges.  From the inequality
$$4n \leq (n+2) + (n+1) + n$$ we see that $n\leq 3$.  But if $n=3$
then $v$ has at most $2(n/2) = n$ positive edges, at most $n-1$
negative edges, and at most $n$ boundary edges, which would lead to the
contradiction that $4n \leq n + (n-1) + n$.

In case (ii), any vertex $v$ of $G_{\a}$ has at most $n+2$ positive
edges, the same number of negative edges, and at most $n$ boundary
edges; so from $4n \leq (n+2) + (n+2) + n$ we see that $n=4$, $|E| =
3$ for all positive interior edges of $\hat G_{\a}$, and $|E| = 4$ for
all boundary edges of $G_{\a}$.  Each label appears three times among
the interior edges endpoints at any vertex $v$ of $G_{\a}$, but since
each of the two families of positive edges at $v$ contains a
Scharlemann cycle, which must all have the same label pair $(1,2)$,
it follows that these labels appear only once among the negative edge
endpoints at $v$, so the label $3$ appears twice among the negative
edge endpoints at $v$.  Since this is true for all vertices $v$ in
$G_{\a}$, it means that the vertex $v_3$ on $G_{\b}$ has $2n_{\a}$
positive edge endpoints, and $n_{\a}$ negative ones, a contradiction.
\endproof

\proclaim{Proposition 5.8} One of the graphs $G_{\a}$ has at most two
vertices.  \endproclaim

\proof Assume $n_1, n_2 \geq 3$.  By Proposition 4.7, $G_{\a}$ is not
special, by Lemma 5.3 $\hat G_{\a}$ does not have an interior vertex
of valency at most 5, and by Lemma 5.4 it cannot have a boundary
vertex of valency at most 4 with a single boundary edge.  Thus by
Proposition 3.1 both $\hat G_1, \hat G_2$ are triangular, which
contradicts Lemma 5.7.  \endproof

\head  \S 6. Nonspecial graphs with $n_1 = 2$ and $n_2 > 2$
\endhead

Throughout this section we will assume that $G_1, G_2$ are not special
graphs.  We will show that the case $n_1=2$ and $n_2 = n > 2$ does not
happen.  Together with Propositions 4.7 and 5.8, this shows that
$n_{\a}$ must be at most 2 for both $\a = 1$ and $2$.

\proclaim{Lemma 6.1} If $n_1 = 2$ then $\hat G_1$ is a subgraph of
that shown in Figure 6.1.  \endproclaim

\proof Since $\hat G_1$ is not a special graph, one of the vertices
$u_1, u_2$ has at most one boundary edge.  If either $u_1$ or $u_2$
does not have a loop, then one can find a vertex $u$ of valency at
most 3 in $\hat G_1$, with at most one boundary edge.  Since each
interior edge represents at most $n$ edges, $u$ would have at least
$2n$ boundary edges, which would imply that $G_2$ is a special graph,
a contradiction.  Hence each vertex $u_i$ has a loop.  It is now easy
to see that $\hat G_1$ must be a subgraph of that in Figure 6.1.
\endproof

\bigskip
\leavevmode

\centerline{\epsfbox{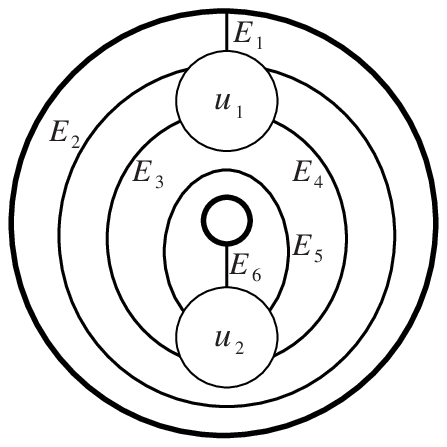}} 
\bigskip
\centerline{Figure 6.1}
\bigskip

Label the edges of $\hat G_1$ as in Figure 6.1.  Denote by $m$ the
number of non-loop interior edges of $G_1$, i.e.\ $m = |E_3| + |E_4|$.

\proclaim{Lemma 6.2}  Suppose $n_1 = 2$, and $n > 2$.

(1) Either $m = 2n$, or $m = 2n -2$ and $E_2$
contains a Scharlemann cycle.  

(2) The two vertices of $G_1$ are antiparallel.  
\endproclaim

\proof (1) If no label appears twice among the endpoints of edges in
$E_2$, then from the labeling on $\bdd u_1$ one can see that either $m
\geq 2n$ or $|E_1| \geq 2n$.  But the second possibility does not
occur because then by Lemma 2.6(2) the graph $G_2$ would be special.
Hence in this case we have $m\geq 2n$.  Since each of $E_3, E_4$
represents at most $n$ edges, we conclude that $m = 2n$.

Now assume that some label appears twice among the endpoints of edges
in $E_2$.  Then $E_2$ contains a Scharlemann cycle $e_1, e_2$, with
label pair $(1,2)$, say.  Since $n > 2$, $E_2$ contains no extended
Scharlemann cycle (Lemma 2.1(5)), so one of these two edges, say
$e_1$, must be an outermost edge among those in $E_2$.  Thus the
endpoints of $e_1$ are either adjacent to those in $E_3\cup E_4$ or to
those in $E_1$.  In the first case, the label sequence of $E_3\cup
E_4$ at $u_1$ is $3, 4, ..., n$, so $m \equiv n-2$ mod $n$.  If $m =
2n-2$ then we are done.  If $m\neq 2n-2$, then since $|E_3|, |E_4|
\leq n$, we must have $m = n-2$.  Thus $|E_1| = \Delta n - m - 2|E_2|
\geq 2n$, which by Lemma 2.6(2) would imply that $G_2$ is special, a
contradiction.  Therefore $e_1$ must be adjacent to $E_1$.  As above,
we have either $|E_1| = 2n-2$, or $|E_1| = n-2$ and $m \geq 2n$.  In
the second case we have $m = 2n$ because $|E_3|, |E_4| \leq n$.  It
remains to show that $|E_1| = 2n-2$ is impossible.

Assume $|E_1| = 2n-2$.  Notice that this happens only if $E_2$
contains a Scharlemann cycle.  Moreover, if $(1,2)$ is the label pair
of the Scharlemann cycle then all labels other than $1,2$ would appear
twice among endpoints of edges in $E_1$.  Thus on $G_2$ each vertex
other than $v_1, v_2$ would have two boundary $1$-edges.  But since
the edges in the Scharlemann cycle and the vertices $v_1, v_2$ form an
essential subgraph of $G_2$, these two parallel $1$-edges must go to
the same boundary component of $A_2$.  By looking at an outermost
vertex one can see that there is a vertex $v_i$ with $i \neq 1,2$, at
which the two boundary $1$-edges are parallel, so they are parallel on
both graphs, contradicting Lemma 2.1(2).  

(2) If $u_1, u_2$ are parallel then $2n-2 \leq |E_3| + |E_4| \leq
2(n/2+1)$, implying that $n=4$ and $|E_3| = |E_4| = 3$.  In this case
both $E_3, E_4$ contain Scharlemann cycles, and by Lemma 2.2(1) they
must have the same label pair $(1,2)$ as the one in $E_2$.  But since
each of the labels $1,2$ appears only once among the endpoints at
$u_1$ of edges in $E_3\cup E_4$, this is impossible.  \endproof

\proclaim{Lemma 6.3} Suppose $n_1 = 2$, and $n > 2$.  Then $G_1$
cannot have $2n$ negative edges.  \endproclaim

\proof We must have $|E_2| > 0$, otherwise $|E_1| \geq 2n$, so $G_2$
would be special, contradicting our assumption.  Assume that $G_1$ has
$2n$ negative edges.  Then $|E_3| = |E_4| = n$, and by Lemma 2.6(3) we
have $|E_1| \leq n$, hence $|E_2| \geq n/2$.  On the other hand, by
Lemma 2.2(2) $|E_2|\leq n/2+1$.  Hence $E_2$ contains either $n/2$ or
$n/2+1$ edges.  We want to show that $|E_2| \neq n/2+1$.  Assuming
otherwise, then since $E_3$ contains $n$ parallel negative edges, by
Lemma 5.1 the graph $G_1$ has at most $n/2$ parallel boundary edges.
On the other hand, we have $|E_1| = \Delta n - m - 2|E_2| \geq n - 2$,
and since $E_2$ contains a Scharlemann cycle with label pair
$(1,2)$, say, $n$ is even.  Therefore we must have $n=4$.  Now in
this case the labels of the edges in $E_1$ are $3,4$, contradicting
Lemma 5.1(2).

\bigskip
\leavevmode

\centerline{\epsfbox{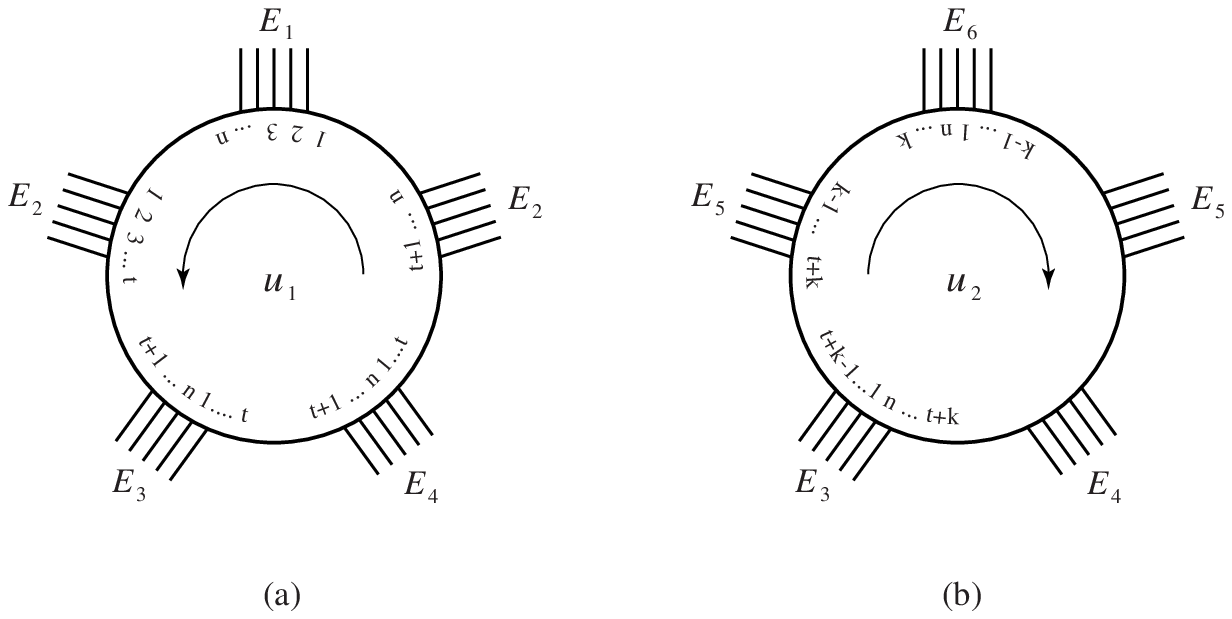}} 
\bigskip
\centerline{Figure 6.2}
\bigskip

We have shown that $|E_2| = n/2$, and $|E_1| = n$.  For the same
reason, we have $|E_5| = n/2$ and $|E_6| = n$.  Without loss of
generality we may assume that $u_1$ is a positive vertex, $u_2$ is
negative, and the edges of $E_1$ have label sequence $1, ..., n$ at
$u_1$.  See Figure 6.2.  Let $t = n/2$.  Since $|E_2| = n/2$, the
label sequence of the endpoints at $u_1$ of the edges of $E_3$ is
$t+1, ..., n, 1, ..., t$.  There is a number $k$ such that the label
sequence at the other end of $E_3$ is $t+k, t+k+1, ..., t+k -1$.  The
number $k \neq 1$, otherwise these edges would be loops in $G_2$, so
$n>2$ would imply that some vertex of $G_2$ does not have a boundary
edge, contradicting the fact that $|E_1| = n$.  Now from Figure 6.2 we
can see that the label sequence of $E_6$ is $k, ..., n, 1, ..., k-1$,
hence the two edges $e_3, e_4$ in $E_6$ labeled $n$ and $1$
respectively, are adjacent (because $k\neq 1$).  Let $e_1, e_2$ be the
edges of $E_1$ labeled $1$ and $n$ respectively.  By Lemma 2.6(3),
each vertex of $G_2$ has at most one family of parallel boundary
edges, so $e_2$ is parallel to $e_3$, and $e_4$ parallel to $e_1$ in
$G_2$.  Let $B(e_1, e_2)$ be the band on $F_1$ between $e_1$ and
$e_2$, and let $B(e_3, e_4)$ be that between $e_3$ and $e_4$.
Similarly, let $B(e_2, e_3)$ and $B(e_4, e_1)$ be the bands on $F_2$
between $e_2, e_3$ and $e_4, e_1$, respectively.  Now we can form an
annulus $A = B(e_1,e_2) \cup B(e_2,e_3) \cup B(e_3,e_4) \cup
B(e_4,e_1)$ in the manifold $M$.  Since the boundary curve $C$ of $A$
on $T_0$ intersects the circle $\bdd v_2$ transversely at a single
point (on the arc $B(e_1, e_2) \cap \bdd u_1$), it is an essential
curve.  This contradicts the fact that the manifold $M$ is
$\bdd$-irreducible and anannular.
\endproof

\proclaim{Lemma 6.4} Suppose $n_1 = 2$, and $n > 2$.  Then $G_1$
cannot have exactly $2n -2$ negative edges.  \endproclaim

\proof If $G_1$ has $2n - 2$ negative edges, then (up to symmetry)
either $|E_3| = |E_4| = n-1$, or $|E_3| = n$ and $|E_4| = n-2$.
Looking at the labeling, one can see that the two loops of $E_2$ near
$E_3\cup E_4$ form a Scharlemann cycle, with label pair $(1,2)$, say.
If $|E_3| =n$ then by Lemma 2.6(3) we have $|E_1| \leq n$, hence
$|E_2| \geq n/2+1$.  Now by Lemma 5.1(1) the graph $G_2$ has at
most $n/2$ boundary vertices, which contradicts the fact that $|E_1| =
n$.  Therefore we must have $|E_3| = |E_4| = n-1$.

For the same reason, the two loops in $E_5$ near $E_3\cup E_4$ form a
Scharlemann cycle, which by Lemma 2.2(1) must have the same label pair
$(1,2)$.  Now we can see that $E_3$ has label sequence $3, 4, ...,
n, 1, $ at $u_1$, and has label sequence $2, 3, ..., n$ at $u_2$.
However, in this case $E_3$ has only one orbit, containing all the
labels, so all the vertices of $G_2$ are parallel to each other, hence
all edges of $G_1$ are negative.  But since $G_1$ contains some
loops, this is a contradiction.  \endproof

\proclaim{Proposition 6.5}  If
$M(r_1), M(r_2)$ are annular, and $\Delta \geq 4$, then $n_{\a} \leq
2$ for $\a = 1,2$.  
\endproclaim

\proof By Proposition 4.7 this is true if one of the $G_{\a}$ is
special.  By Proposition 5.8 one of the graphs, say $G_1$, has at most
two vertices.  Since the two possibilities in Lemma 6.2(1) have been
ruled out by Lemmas 6.3 and 6.4, the case $n_1\leq 2$ and $n_2 > 2$
cannot happen.  \endproof

\head   \S 7.  Special graphs with $n_1=1$ and $n_2 = 2$
\endhead

\proclaim{Proposition 7.1} If $G_{\a}$ is special, then $\Delta = 4$,
up to relabeling $n_1 = 1$, $n_2 = 2$, and the manifold $M$ is the
exterior of the Whitehead link.  \endproclaim

\proof By Lemma 4.1, both graphs must be special.  By Proposition 4.7,
up to relabeling we must have $n_1=1$, $n_2 = 2$, and $G_1$ has
exactly two interior edges $e_1,e_2$.

Assume $\Delta = 5$.  By Lemma 4.4 the jumping number $q=1$.  There is
a pair of adjacent boundary $1$-edges $e_1, e_2$ at $v_1$ in $G_2$,
which by Lemma 2.5(1) should also be adjacent at $u_1$ in $G_1$ among
all $1$-edges; but since the two families of boundary edges at $u_1$
are separated by two interior edges, $e_1, e_2$ must be in the same
family, so they are parallel on both graphs, a contradiction.
Therefore we must have $\Delta = 4$.

Now the Whitehead link exterior $W$ does admit two annular Dehn
fillings $W(r_1), W(r_2)$ with $\Delta(r_1,r_2) =4$, $n_1 = 1$, and
$n_2 = 2$, see [GW1, Theorem 7.5].  It remains to show that the
manifold satisfying these conditions is unique.

Each vertex of $G_2$ has two boundary edges, which are nonparallel
because $G_2$ is special.  Thus the graph $G_2$ must be as shown in
Figure 7.1(b).  Similarly, since $G_1$ is special it has two
families of parallel boundary edges.  The loops have different labels
at their two endpoints, so each family of boundary edges of $G_1$
contains an even number of edges.  Hence $G_1$ must be as shown in
Figure 7.1(a).

\bigskip
\leavevmode

\centerline{\epsfbox{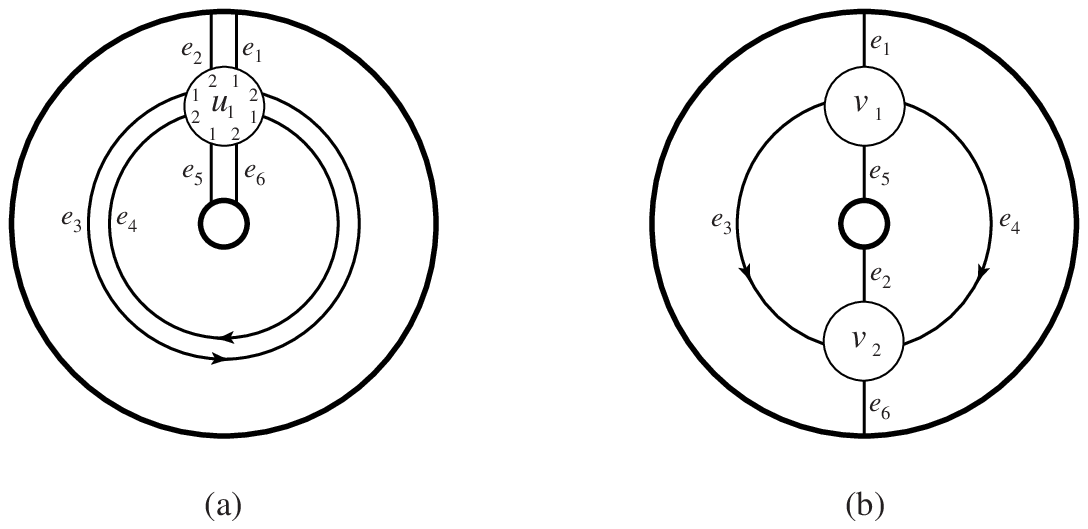}} 
\bigskip
\centerline{Figure 7.1}
\bigskip

Label the six edges of $G_1$ as in the figure.  Orient $e_3, e_4$ so
that on $G_1$ they have label $1$ at their tails.  Up to symmetry we
may assume that the edge $e_1$ on $G_2$ is as shown in Figure 7.1(b).
The label $1$ endpoints of edges $e_1, e_3, e_5, e_4$ appear
successively on $\bdd u_1$, hence by Lemma 2.5(1) they also appear in
this order on $\bdd v_1$ in $G_2$, so these edges must be as shown in
the figure.  Similarly by looking at the label $2$ endpoints of $e_2,
e_4, e_6, e_3$ one can determine the edges $e_4$ and $e_6$.  Therefore
up to symmetry the graphs $G_{\a}$ are exactly as shown in the figure.
We need to show that these graphs uniquely determine the manifold $M$.

Recall that $F_{\a}$ denotes the punctured annulus $A_{\a} \cap M$.
Let $X = N(F_1 \cup T_0)$, and let $Y = N(F_1\cup F_2 \cup T_0)$,
where the regular neighborhoods are taken in $M$.  The frontier of $X$
in $M$, i.e.\ $X \cap \overline{M-X}$, is a surface $F$, which is a
four punctured sphere.  Note that $Y$ is obtained from $X$ by adding
regular neighborhoods of the faces of $G_2$.  Each of the four faces
of $G_2$ is a disk $D_i$ with $\bdd D_i = c_i \cup c_i'$, where $c'_i$
is an arc on $\bdd M$, and $c_i$ an arc on $F_1\cup T_0$.  Let $\tilde
c_i$ be the arc $D_i \cap F$.  Then the frontier of $Y = X \cup (\cup
N(D_i))$ in $M$ is a properly embedded surface $F'$, homeomorphic to
the surface obtained by cutting $F$ along the arcs $\tilde c_i$.  Thus
$Y$ and $X$ are homeomorphic, but they are embedded in $M$
differently.  Note that $Y$ is uniquely determined by the graphs
$G_1$ and $G_2$.

It is easy to see that all the $\tilde c_i$ are essential arcs on $F$.
Since each boundary component of $F$ meets $\cup \tilde c_i$ twice,
after cutting along all these $\tilde c_i$, the remnant, and hence
$F'$, consists of either two disks, or two disks and an annulus.  In
fact, by examining the graphs, one can see that $F'$ indeed consists
of two disks and an annulus.  Since $M$ is irreducible and
$\bdd$-irreducible, the disk components of $F'$ are boundary parallel.
If the annular component $A$ of $F'$ is incompressible in $M$ then $A$
is also boundary parallel because $M$ is anannular and irreducible, so
$M$ would be homeomorphic to $Y$, which in turn is homeomorphic to
$X$.  Let $C$ be an essential curve on $T_0$ disjoint from $\bdd F_1$.
Then $C \times I$ in $T_0 \times I$ would be an essential annulus in
$X$, contradicting the fact that $M$ is anannular.  Therefore $A$ must
be compressible.  Let $D$ be a compressing disk of $A$ in $M$.  Then
$D$ lies in either $Y$ or $M- \Int Y$.  We show that the first case is
impossible.

First notice that the surface $F_1$ cuts $X$ into a manifold $F \times
I$, in which both $F$ and the two copies of $F_1$ are incompressible.
By an innermost circle argument one can show that $F$ is
incompressible in $X$.  Under the homeomorphism $Y \cong X$, $A$ can
be considered as a subsurface of $F$, hence $A$ is also incompressible
unless the core of $A$ is a trivial curve on $F$.  On the other hand,
notice that $A$ is a component of $\bdd Y - \Int F''$, where $F'' = Y
\cap (\bdd M - T_0)$ is a neighborhood of $\bdd A_1 \cup \bdd A_2$ on
$\bdd M$, which is connected.  Since $\bdd A \subset F''$, it follows
that the core of $A$ is nonseparating on $\bdd Y$, hence it is
nontrivial on $F$.  This completes the proof that $A$ is
incompressible in $Y$.

Hence the compressing disk $D$ of $A$ lies in $M-\Int Y$.  Let $M'$ be
the union of $Y$ and a regular neighborhood of $D$.  Then the frontier
of $M'$ in $M$ is a set of disks, which must be boundary parallel
because $M$ is irreducible and $\bdd$-irreducible.  Therefore $M'$ is
homeomorphic to $M$.  It follows that $M$ is obtained from $Y$ by
adding a 2-handle along the core of $A$, and hence is uniquely
determined by the graphs $G_1$ and $G_2$.  \endproof

\head  \S 8.  Nonspecial graphs with $n_{\a} \leq 2$
\endhead

First note that if $n_{\a} = 1$ and $G_{\a}$ is not special, then the
unique vertex of $G_{\a}$ has valency at most $3$ in $\hat G_{\a}$,
and hence by Lemma 2.2(3) $G_{\a}$ has at least $2n_{\b}$ parallel
boundary edges.  By Lemma 2.6(2) this implies that $G_{\b}$, and
therefore (by Lemma 4.1) $G_{\a}$, is special, a contradiction.  Hence
if $G_1, G_2$ are not special and $n_1, n_2 \leq 2$, we must have $n_1
= n_2 = 2$.

\proclaim{Lemma 8.1} Suppose that $n_1 = n_2 = 2$ and $G_1, G_2$ are
not special.  Then for $\a = 1,2$, the two vertices of $G_{\a}$ are
antiparallel, $\hat G_{\a}$ is a subgraph of the graph $\hat G$ in
Figure 6.1, and one of the following holds.

(i) $\Delta = 4$, each interior edge of $\hat G_{\a}$ represents two
edges of $G_{\a}$, and $G_{\a}$ has no boundary edges.

(ii) $\Delta = 5$, each edge of $\hat G_{\a}$ represents two edges of
$G_{\a}$, and the jumping number $q=2$.  \endproclaim

\proof By Lemma 6.1, $G_{\a}$ is a subgraph of the graph $\hat G$
shown in Figure 6.1.  Each vertex $v$ of $G_{\a}$ must have a loop,
otherwise some vertex would have valency $3$ in $\hat G_{\a}$ with a
single boundary edge, so by Lemma 2.6(2) $G_{\b}$ would be special,
contradicting the assumption.
Since a loop in $G_{\a}$ is a non-loop negative edge of $G_{\b}$, it
follows that each graph $G_{\b}$ has some negative edges, hence the
two vertices of $G_{\b}$ must be antiparallel, $\b=1,2$.  By Lemma
2.2(3) each interior edge of $\hat G_{\a}$ represents at most two
edges of $G_{\a}$.  Similarly, each boundary edge of $\hat G_{\a}$
also represents at most two edges of $G_{\a}$, by Lemma 2.1(2).

First assume $\Delta = 4$.  Notice that a vertex of $G_{\a}$ has
either no boundary edge or two boundary edges, for if it has exactly
one boundary edge then the loops based at that vertex would have the
same label at their two endpoints, which contradicts the parity rule.
Since two boundary edges at a vertex of $G_{\a}$ correspond to
boundary edges at different vertices of $G_{\b}$, it follows that
either both vertices of $G_{\a}$ have two boundary edges, or they both
have no boundary edges.  The second possibility gives rise to
conclusion (i) in the lemma.

Assume that each vertex of $G_{\a}$ has two boundary edges.  Then
there are a total of 6 interior edges in each graph.  Note that an
interior edge is a loop on $G_{\a}$ if and only if it is a non-loop on
$G_{\b}$ because of the parity rule, hence one of the graphs, say
$G_1$, has at least three loops.  Without loss of generality we may
assume that there are two loops $e_1, e_2$ based at the vertex $u_1$.
Consider their label $1$ endpoints.  Because there are two boundary
edges at $u_1$, these two endpoints are non adjacent among all label
$1$ endpoints at $u_1$.  Now look at the graph $G_2$.  By Lemma 2.5(1)
$e_1, e_2$ are non adjacent $1$-edges at $v_1$ among all $1$-edges.
However, since they are non-loops in $G_2$, they are contained in the
two adjacent families $E_3, E_4$ in Figure 6.1.  Since $E_3\cup E_4$
contains a total of at most four edges, $e_1, e_2$ are adjacent among
all $1$-edges at $v_1$.  This contradiction completes the proof of the
lemma for the case $\Delta = 4$.

Now assume $\Delta = 5$.  Since each vertex of $\hat G$ has valency 5,
and since each edge of $\hat G_{\a}$ represents at most two edges of
$G_{\a}$, $\Delta = 5$ implies that each edge of $\hat G_{\a}$
represents exactly two edges.  By the same argument as above one can
show that the jumping number $q$ cannot be $1$, so we are in case (ii).
\endproof

\proclaim{Lemma 8.2}  There is a unique irreducible, 
$\bdd$-irreducible, anannular manifold $M$ which
admits two annular Dehn fillings $M(r_1), M(r_2)$ with $\Delta(r_1,
r_2) = 5$.
\endproclaim

\proof By Lemma 8.1, the graphs must be as shown in Figure 8.1.
We first show that the edge correspondence and the labelings of the 
vertices are unique up to symmetry.

Reflecting the annuli vertically and changing their orientations if
necessary, we may assume that the vertices $u_1, v_1$ are positive,
and the labeling of edge endpoints at $\bdd u_1, \bdd v_1$ are as
shown.  Any non-loop edge has the same label on its two endpoints,
because it is a loop edge on the other graph.  Thus the labeling on
$\bdd u_2, \bdd v_2$ is determined by that on $\bdd u_1, \bdd v_1$,
respectively.  Orient the edges so that a non-loop edge goes from
$u_1$ to $u_2$ (resp.\ $v_1$ to $v_2$).  Then dually the orientation
of a loop edge must go from label $1$ to label $2$.  Label the edges
of $G_1$ as in Figure 8.1(a).

If $P_1, ..., P_5$ are the points of $u_1\cap v_1$, appearing in this
order on $\bdd u_1$ along its orientation, then since the jumping
number $q=2$, they appear in the order $P_1, P_3, P_5, P_2, P_4$ on
$\bdd v_1$ either along or against the orientation of $\bdd v_1$.  In
other words, along the orientation of $\bdd v_1$ they either appear in
this order, or in the order $P_1, P_4, P_2, P_5, P_3$.  In the second
case, write $(Q_1, Q_2, ..., Q_5 ) = (P_1, P_4, P_2, P_5, P_3)$; then
$(P_1, ..., P_5) = (Q_1, Q_3, Q_5, Q_2, Q_4)$.  Hence by interchanging
the roles of $G_1$ and $G_2$ if necessary, we may assume that the
points appear as $(P_1, P_3, P_5, P_2, P_4)$ on $\bdd v_1$ along the
orientation of $\bdd v_1$.

Now we can see that the labeling of the edges on $G_2$ is completely
determined by that of $G_1$: The $1$-edges at $u_1$ appear in the
order $a, c, e, k, d$ in the positive direction, so at $v_1$ they appear
in the order $a, e, d, c, k$, where $a$ is the unique boundary edge at
$v_1$ labeled $1$.  The order of the 2-edges at $u_1$ is $b,d,f,l,c$,
so dually the $1$-edges at $v_2$ are in the order $b, f, c, d, l$.
Similarly by looking at $u_2$ one can determine the labeling of the
remaining edges in $G_2$.  See Figure 8.1(b).

\bigskip
\leavevmode

\centerline{\epsfbox{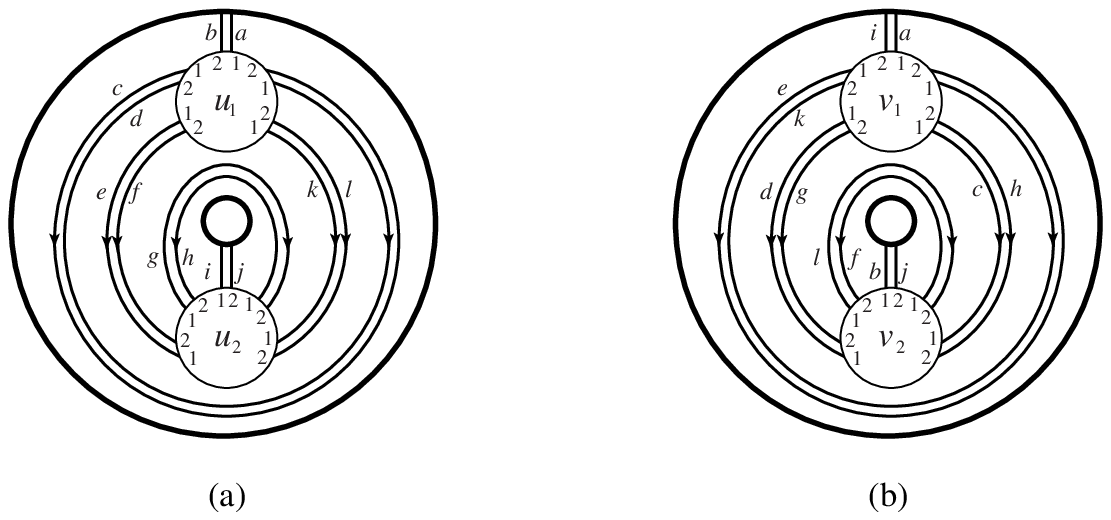}} 
\bigskip
\centerline{Figure 8.1}
\bigskip

It remains to show that the manifold $M$ is uniquely determined by
these graphs.  As in the proof of Proposition 7.1, consider the
submanifold $X = N(A_1 \cup J_1)$ of $M(r_1)$.  Since $J_1$ intersects
$A_1$ in two meridian disks of opposite sign, the frontier $F$ of $X$
consists of two components $F_b, F_w$, each being a twice punctured
torus, called the {\it black surface\/} and the {\it white surface\/}
respectively.  A face of $G_2$ is black or white according to whether
it intersects the black surface or the white surface.  Note that each
face of $G_2$ intersects $F$ in a circle or an arc, so it is either
black or white, but not both.

Let $D_1$ be a face of $G_2$ bounded by a pair of parallel loops, and
let $D_2$ be the triangular interior face of $G_2$ adjacent to $D_1$.
Since they have an edge in common, they are of different colors, so we
may assume that $D_1$ is black and $D_2$ is white.  The boundary of
$D_1$ intersects a meridian of $J_1$ twice in the same direction,
hence $\bdd D_1$ is a nonseparating curve on $F_b$.  After adding a
neighborhood of $D_1$ to $X$, the black frontier is homeomorphic to
the surface obtained by 2-surgery on $F_b$ along $\bdd D_1$, hence is
an annulus $A_b$.  Since its boundary components are essential curves
on $\bdd M$, and since $M$ is $\bdd$-irreducible, $A_b$ is
incompressible in $M$, and hence is boundary parallel in $M$.
Similarly, since the boundary of $D_2$ intersects a meridian of $J_1$
three times, $\bdd D_2$ is a nonseparating curve on $F_w$, so after
adding $N(D_2)$ the white frontier becomes an annulus $A_w$, which for
the same reason must be boundary parallel in $M$.  It follows that $M$
is homeomorphic to $N(F_1 \cup T_0 \cup D_1 \cup D_2)$, where $F_1$ is
the punctured annulus $A_1 \cap M$.  The boundary curves of $D_i$ are
determined by the graphs, which have been determined (up to symmetry)
as above.  Hence the manifold $M$ is uniquely determined.  \endproof

\proclaim{Lemma 8.3} There is a unique irreducible,
$\bdd$-irreducible, anannular manifold $M$ which admits two annular
Dehn fillings $M(r_1), M(r_2)$ with $\Delta(r_1, r_2) = 4$ and $n_1 =
n_2 = 2$.  \endproclaim

\proof The proof is similar to that of Lemma 8.2. In this case the
jumping number is $1$, and one can show that up to symmetry the graphs
must be as shown in Figure 8.2.  The proof that $M$ is determined by
the graphs is the same as in the proof of Lemma 8.2.
\endproof

\bigskip
\leavevmode

\centerline{\epsfbox{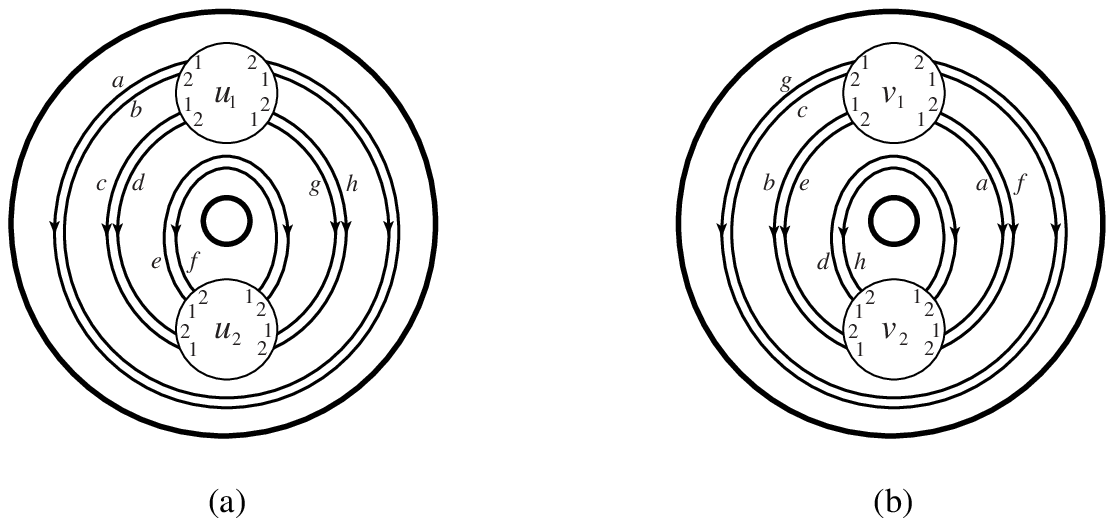}} 
\bigskip
\centerline{Figure 8.2}
\bigskip

We now prove Theorem 1.1, which we restated here for the reader's
convenience.  

\proclaim{Theorem 1.1} Suppose $M$ is a compact, connected,
orientable, irreducible, $\bdd$-irreducible, anannular $3$-manifold
which admits two annular Dehn fillings $M(r_1)$, $M(r_2)$ with $\Delta
= \Delta(r_1, r_2) \geq 4$.  Then one of the following holds.

(1)  $M$ is the exterior of the Whitehead link, and $\Delta = 4$.

(2)  $M$ is the exterior of the $2$-bridge link associated to the
rational number $3/10$, and $\Delta = 4$.

(3)  $M$ is the exterior of the $(-2,3,8)$ pretzel link, and $\Delta =
5$.  
\endproclaim

\proof
By Proposition 6.5, we must have $n_{\a} \leq 2$ for $\a = 1,2$.  If
$G_{\a}$ is special, then by Proposition 7.1 the manifold $M$ is the
exterior of the Whitehead link.  If $G_{\a}$ is nonspecial, then by
Lemma 8.1 the graphs $G_{\a}$ must be as in Figure 8.1 or 8.2, and by
Lemmas 8.2 and 8.3, in each case the manifold $M$ is uniquely
determined by the graphs; hence there are at most three manifolds $M$
which may admit two annular Dehn fillings of distance at least $4$
apart.  On the other hand, it has been shown in [GW1, Theorem 7.5]
that each of these manifolds admits two such fillings.  Hence the
result follows.
\endproof

\enddocument